\documentclass[12pt]{article}
\usepackage{a4wide,latexsym,amsmath,amssymb,amsthm}

\newtheorem{theorem}{Theorem}

\newtheorem{lemma}{Lemma}

\begin{document}

\title{\vspace{-3cm} Chiral covers of regular maps of given type}
\author{}
\date{}
\maketitle

\begin{center}
\vspace{-1.3cm}

{\large Olivia Reade} \\
\vspace{1mm} {\small Open University, Milton Keynes, UK}\\

\vspace{4mm}

{\large Jozef \v Sir\'a\v n} \\
\vspace{1mm} {\small Slovak University of Technology, Bratislava, Slovakia}

\vspace{4mm}

\end{center}

\begin{abstract}
With the help of the theory of holomorphic and anti-holomorphic differentials, G. A. Jones [Chiral covers of hypermaps, Ars Math. Contemp. 8 (2015), 425--431] proved that every regular hypermap of a non-spherical type is covered by an infinite number of orientably-regular but chiral hypermaps of the same type. We present a different proof of the same result for regular {\em maps}, based on parallel products of maps and existence of chiral maps of a given hyperbolic type with a symmetric or an alternating automorphism group.
\end{abstract}

\section{Introduction}\label{sec:intro}

An {\em orientably-regular map} $M$ is a cellular embedding of a connected graph on an orientable surface, such that the group $G={\rm Aut}^+(M)$ of all orientation-preserving automorphisms of the embedding is transitive (and hence regular) on the mutually incident vertex-edge pairs. In such a map, all face boundaries have the same length, say, $m$, and all vertex valencies are also the same, say, $n$; we then speak of a map of {\em type} $\{m,n\}$. The type, and the map, are said to be {\em non-spherical} if $1/m+1/n\le 1/2$, and {\em hyperbolic} if $1/m+1/n < 1/2$. If $M$ admits also an automorphism that reverses the orientation of its carrier surface, then the map is {\em reflexible}; otherwise it is {\em chiral}.
\smallskip

Given a pair of orientably-regular maps $M$ and $\tilde M$ of the same type, the map $\tilde M$ is a {\em smooth cover} of $M$ if there exists a smooth (that is, unbranched) covering of the carrier surface of $M$ by the carrier surface of $\tilde M$ that preserves the cell structure of both maps.
\smallskip

In 2015, G. A. Jones  proved a remarkable result \cite[Theorem 2.1]{J1} which implies that every orientably-regular reflexible map of an arbitrary non-spherical type has an infinite number of orientably-regular but {\em chiral} smooth covers. Theorem 2.1 of \cite{J1} is actually more general and concerns {\em hypermaps}, which will be mentioned later. Even more remarkably, one of the tools used in the proof of this result is the theory of holomorphic and antiholomorphic differentials on Riemann surfaces.
\smallskip

In Section \ref{sec:main} of this note we offer a different proof of the theorem of Jones for regular {\em maps},  combining a method of building new orientably-regular maps from old, known as the parallel product of maps, with existence of orientably-regular but chiral maps of arbitrary non-spherical type having symmetric or alternating automorphism groups. In the process it was necessary to consider a number of fine points concerning subgroups of index $2$ in triangle groups (Section \ref{sec:special}) and relating direct product of groups and their $2$-extensions to their automorphism groups (Section \ref{sec:para}); this all is preceded by a brief introduction to the theory of orientably-regular maps in Section \ref{sec:OrReg}.

\section{Orientably-regular maps: basic facts revisited}\label{sec:OrReg}

In this section we recall a few basic facts about orientably-regular maps and various ways of their representation. We will also discuss chirality and reflexibility, followed by covers and parallel products of such maps; for details we refer the reader to the fundamental paper \cite{JoSi} and the survey \cite{Si-surv}.
\smallskip

Let $M$ be a finite orientably-regular map of type $\{m,n\}$. Choose in $M$ a vertex-edge incident pair $(v,e)$ and a face $F$ containing both $v$ and $e$. By orientable regularity of $M$, there is an automorphism $x\in G={\rm Aut}^+(M)$ acting locally as an $m$-fold rotation of $M$ about the centre of $F$, together with an automorphism $y\in G$ acting locally as an $n$-fold rotation of $M$ about the vertex $v$, such that their composition $xy$ acts locally as a $2$-fold rotation of $M$ about the centre of the edge $e$. Because of these `rotary' automorphisms, such maps are also known as {\em rotary}, cf. \cite{Wil2}.   
\smallskip

Connectedness of the underlying graph of $M$ implies that the elements $x$ and $y$ generate the group $G$, so that it admits a presentation of the form
\begin{equation}\label{eq:Aut+}
G={\rm Aut}^+(M) = \langle x,\, y\ |\ x^m,\,y^n,\, (xy)^2,\ldots\rangle
\end{equation}
It follows that $G$ is a smooth quotient of the {\em $(m,n,2)$-triangle group} $\Delta(m,n,2)$, given in the form of a $2$-generator and $3$-relator presentation
\begin{equation}\label{eq:Delta}
\Delta(m,n,2) = \langle X,\,Y\ |\ X^m,\,Y^n,\,(XY)^2\rangle
\end{equation}
The kernel of the natural smooth group epimorphism $\Delta(m,n,2)\to G$ is known as the {\em map subgroup} of $M$; it is a normal, torsion-free subgroup of $\Delta(m,n,2)$ of finite index.
\smallskip

By the general theory of orientably-regular maps \cite{JoSi}, there is a one-to-one correspondence between (isomorphism classes of) orientably-regular maps of type $\{m,n\}$, torsion-free normal subgroups of $\Delta(m,n,2)$ of finite index, and presentations of $2$-generated finite groups of the form \eqref{eq:Aut+}. An orientably-regular map $M$ of type $\{m,n\}$ with $G={\rm Aut}^+(M)$ as in \eqref{eq:Aut+} will be denoted $M={\rm Map}(G;x,y)$.
\smallskip

The well-known construction of a dual map to a given map has the following algebraic equivalent: if $M={\rm Map}(G;x,y)$ is an orientably-regular map of type $\{m,n\}$, its {\em dual map} $\tilde{M}={\rm Map}(G;x',y')$, of type $\{n,m\}$, is obtained from the presentation \eqref{eq:Aut+} of the {\em same group} $G$ by letting $x'=y$ and $y'=x$. Informally, duality in orientably-regular maps is induced by interchanging the roles of the rotary automorphisms $x$ and $y$, or, at the level of triangle groups, by swapping the roles of $X$ and $Y$.
\smallskip

Reflexibility, that is, existence of a map automorphism reversing the orientation of the carrier surface, has also a well-known equivalent description in our context: An orientably-regular map $M={\rm Map}(G;x,y)$ is reflexible if and only if there is a {\em group automorphism} of $G$ that inverts both $x$ and $y$. It follows that in order to prove that such a map $M$ is chiral it is sufficient to {\em exclude} the existence of an automorphism of $G$ inverting its two generators.
\smallskip

Let $M_1$ and $M_2$ be orientably-regular maps of the same type $\{m,n\}$, given, respectively, by their map subgroups $L_1,\,L_2\triangleleft \Delta(m,n,2)$. We say that $M_2$ is a {\em smooth regular cover} (or, simply, a {\em smooth cover}) of $M_1$ if $L_2$ is a (normal) subgroup of $L_1$. This situation exactly corresponds to a smooth (that is, unbranched) regular topological covering of the carrier surface of $M_1$ by a carrier surface of $M_2$, which, in addition, preserves the cell structure of both maps. If one represents the two maps by group presentations as in \eqref{eq:Aut+}, that is, if one lets $M_i={\rm Map}(G_i;x_i,y_i)$ for $i=1,2$, then $M_2$ is a smooth cover of $M_1$ if and only if there is an epimorphism $G_2\to G_1$ sending $x_2$ to $x_1$ and $y_2$ to $y_1$. This follows from the way the two representations are related: namely, by $G_i\cong \Delta(m,n,2)/L_i$ for $i=1,2$
\smallskip

Finally, we address a useful way of constructing new maps from old, focusing on the special case when the maps are orientably-regular and of the same type. Let $M_1$ and $M_2$ be orientably-regular maps, both of type $\{m,n\}$, with respective map subgroups $L_1,\,L_2$ of $\Delta(m,n,2)$. The {\em parallel product} $M_1||M_2$ of the two maps is the orientably-regular map determined by the map subgroup $L_1\cap L_2$. We note that since both $L_1$ and $L_2$ are torsion-free, normal subgroups of the triangle group $\Delta(m,n,2)$, their intersection has the same properties and hence $M_1||M_2$ is again of type $\{m,n\}$. Note that the map $M_1|| M_2$ is a common smooth cover of both $M_1$ and $M_2$.
\smallskip

An equivalent description of the parallel product of $M_1$ and $M_2$ can be given in terms of representing the two maps in the form $M_i={\rm Map}(G_i;x_i,y_i)$ for $i=1,2$. Letting $G$ denote the subgroup of the direct product $G_1 \times G_2$ generated by the elements $x=(x_1,x_2)$ and $y=(y_1,y_2)$, then $M_1||M_2 = {\rm Map}(G;x,y)$. The group $G$ itself is also called a {\em parallel product} of $G_1$ and $G_2$; it is a special instance of a sub-direct product of groups.

\section{Exceptional orientably-regular maps}\label{sec:special}

An orientably-regular map $M$ with map subgroup $L\triangleleft \Delta =\Delta(m,n,2)$ will be said to be {\em exceptional} if $L$ is also contained in some subgroup, say, $\Gamma$, of index $2$ in $\Delta$; more specifically, we will occasionally say that $M$ is {\em $\Gamma$-exceptional}. Since there will be a need to consider such maps, we will have a brief look at some related facts.
\smallskip

Let $\Gamma$ be a (normal) subgroup of $\Delta=\langle X,Y\rangle$ of index $2$; in symbols, $\Gamma\, \triangleleft_{\,2} \,\Delta$. Then, either $X\in \Gamma$ but $Y\notin \Gamma$ (and then $Y^2\in \Gamma$, implying that $n$ is even), or $X\notin \Gamma$ but $Y\in \Gamma$, (and then $X^2\in \Gamma$, so that $m$ is even) or else $X,Y\notin \Gamma$ (in which case $X^2,Y^2\in \Gamma$, and both $m,n$ are even). These three situations will be designated, respectively, by specifying $\Gamma$ in the form $\Gamma(X,Y^2)$, $\Gamma(X^2,Y)$ and $\Gamma(X^2,Y^2)$. Interchange of the roles of $X$ and $Y$ in \eqref{eq:Delta} induces duality in orientably-regular maps, so that the cases $\Gamma(X,Y^2)$ and $\Gamma(X^2,Y)$ can be treated as dual to each other. We will therefore consider only the subgroups $\Gamma(X,Y^2)$ (if $n$ is even) and $\Gamma(X^2,Y^2)$ (if both $m,n$ are even) in what follows; referring to them as subgroups of type $A$ and $B$, respectively. They are well known and their presentations can easily be derived by the Reidemeister-Schreier method:
\begin{equation}\label{eq:Y2}
{\rm Type\ A:} \ \ \Gamma(X,Y^2) = \langle U,\,V\ |\ U^m,\,V^{n/2},\,(UV)^m\rangle, \ \  {\rm and} \ \ \ \ \ \ \ \ \ \ \
\end{equation}
\begin{equation}\label{eq:X2Y2}
{\rm Type\ B:}\ \ \Gamma(X^2,Y^2) = \langle U,\,V,\,W\ |\ (UW)^{m/2},\,V^{n/2},(UV)^2,\,W^2\rangle,
\end{equation}
where $V=Y^2$ in both cases, while in \eqref{eq:Y2} one has $U=YXY^{-1}$, $UV=X^{-1}$, and in \eqref{eq:X2Y2} the new generators satisfy $U=XY^{-1}$ and $W=YX$.
\smallskip

Let now $M={\rm Map}(G;x,y)$ be a $\Gamma$-exceptional orientably-regular map for a subgroup $L\,\triangleleft \,\Delta = \Delta(m,n,2)$, with at least one of $m$, $n$ even. Up to duality, we may (and will) assume that $n$ is even in what follows. The subgroup $\Gamma\,\triangleleft_{\,2}\,\Delta $ containing $L$ is then of type A or B, that is,  one of $\Gamma(X,Y^2)$, $\Gamma(X^2,Y^2)$. We denote by $H$ the image in $G$ of the restriction $\theta_\Gamma$ of the natural projection $\theta: \Delta \to G$ given by $X\theta=x$ and $Y\theta =y$ to the subgroup $\Gamma$. Clearly, $H\triangleleft_{\,2} G$, and \eqref{eq:Y2} and \eqref{eq:X2Y2} imply that $H$ admits the following presentation:
\begin{equation}\label{eq:H-Y2}
H = \langle u,\,v\ |\ u^m,\,v^{n/2},\,(uv)^m,\ldots \rangle \ \ {\rm if}\ \Gamma=\Gamma(X,Y^2), \ \ {\rm and}
\end{equation}
\begin{equation}\label{eq:H-X2Y2}
H = \langle u,\,v,\,w\ |\ (uw)^{m/2},\,v^{n/2},(uv)^2,\,w^2,\ldots \rangle  \ \ {\rm if}\ \Gamma=\Gamma(X^2,Y^2),
\end{equation}
with $v=y^2$ in both cases, and in \eqref{eq:H-Y2} the generators satisfy $u=yxy^{-1}$, $uv=x^{-1}$, while in \eqref{eq:X2Y2} one has $u=xy^{-1}$ and $w=yx$.
\smallskip

We will now address extensions of certain automorphisms of $H$ to the group $G$.
\smallskip

\begin{lemma}\label{lem:HtoG}
Let $M={\rm Map}(G;x,y)$ be an exceptional orientably-regular map of type $\{m,n\}$ for some even $n$ and for some subgroup  $\Gamma\,\triangleleft_{\,2}\, \Delta(m,n,2)$, of type {\rm A} or {\rm B}. Let $H\, \triangleleft_{\,2}\, G$ be the image of the restriction of the natural epimorphism $\Delta\to G$ to $\Gamma$.
\medskip

\noindent {\rm (a)} If $\Gamma$ is of type {\rm A} and if $\vartheta$ is an automorphism of $H=\langle u,v\rangle$, with $u=yxy^{-1}$ and $v=y^2$, such that $u\vartheta = v^{-1}u^{-1}v$ and $v\vartheta = v^{-1}$, then $\vartheta$ extends to an automorphism of $G$ inverting both $x$ and $y$, so that the map $M$ is reflexible.
\medskip

\noindent {\rm (b)} If $\Gamma$ is of type {\rm B} and if $\vartheta$ is an automorphism of $H=\langle u,v,w\rangle$, with $u=xy^{-1}$, $v=y^2$ and $w=yx$, such that $u\vartheta = wv$, $v\vartheta = v^{-1}$ and $w\vartheta = uv$, then $\vartheta$ extends to an automorphism of $G$ inverting both $x$ and $y$, and hence $M$ is again reflexible.
\end{lemma}

{\bf Proof.} We only prove the part (b); the proof of (a) is similar and even easier. Let $\vartheta$ be an automorphism of $H=\langle u,v,w\rangle$ satisfying $u\vartheta = wv$, $v\vartheta = v^{-1}$ and $w\vartheta = uv$. The groups $G=\langle x,y\rangle$ and $H=\langle u,v,w\rangle\, \triangleleft_{\,2}\, G$ are related by $G=H\langle y\rangle$, with $y^2\in G$; in particular, $y$ conjugates $H$ to itself.
\smallskip

We first show that for every $g\in H$, the relation $h=ygy^{-1}$ implies $h\vartheta = y^{-1}(g\vartheta)y$. Clearly, it is sufficient to verify this in the case when $g$ is any of the three generators $u,v,w$ of $H$, and in part (b) one has $u=xy^{-1}$, $v=y^2$ and $w=yx$. If $g=u=xy^{-1}$, then $h=ygy^{-1}=yxy^{-2}=wv^{-1}$, but at the same time for $g\vartheta = u\vartheta = wv$ one has $y^{-1}(g\vartheta) y = y^{-1}(wv)y = y^{-1}(yxy^2)y = xy^{-1}y^4=uv^2$, which is equal to $h\vartheta = (wv^{-1})\vartheta = (w\vartheta)(v^{-1}\vartheta) = uv^2$. The calculations for $g\in \{v,w\}$ are entirely similar and omitted.
\smallskip

As every element of $G{\setminus}H$ is uniquely expressible in the form $hy$ for $h\in H$, one may try to extend $\vartheta$ to $G{\setminus}H$ by letting $(hy)\vartheta = (h\vartheta)y^{-1}$. As a mapping, $\vartheta$ is a well defined bijection on $G$ inverting $y$, and it remains to show that it also extends to a group homomorphism of $G$. To this end, it is sufficient to verify that this extension of $\vartheta$ preserves the product $yg$ for every $g\in H$, that is, to verify that $(yg)\vartheta = (y\vartheta)(g\vartheta) = y^{-1}(g\vartheta)$. Realising that $yg=hy$ for $h=ygy^{-1}\in H$ and using the observation $h\vartheta = y^{-1}(g\vartheta)y$ from the previous paragraph for such $g$ and $h$, one obtains $(yg)\vartheta = (hy)\vartheta = h\vartheta y^{-1} = y^{-1}(g\vartheta) = (y\vartheta) (g\vartheta)$, as expected.
\smallskip

It remains to verify that the extended automorphism $\vartheta$ of $G$ also inverts $x$. But observing that $x=uy$ one has $x\vartheta = u\vartheta y^{-1} = wvy^{-1} = yxy^2y^{-1}= yxy = x^{-1}$.  \hfill $\Box$

\section{Parallel product of groups: a special case}\label{sec:para}

The key point towards our main result will be relating the automorphism group $G\cong {\rm Aut}^+(M_1||M_2)$ of the parallel product to the automorphism group $G_i\cong {\rm Aut}^+(M_i)$ of the constituents, for $i=1,2$. We do this by a lemma which is likely to be folklore and valid in a more general setting. In the statement, $\Delta$ is assumed to be an arbitrary group, and $\Gamma$ an arbitrary subgroup of $\Delta$ of index $2$; in the next section we will apply the result to triangle groups and their subgroups of index $2$.

\begin{lemma}\label{lem:parallel}
Let $K$ and $L$ be distinct normal subgroups of a group $\Delta$.
\medskip

\noindent {\rm (a)} If $KL=\Delta$, then $\Delta/(K\cap L)\cong \Delta/K\times \Delta/L$.
\medskip

\noindent {\rm (b)} Let $\Gamma$ be a normal subgroup of index $2$ of $\Delta$, such that $\Delta=\Gamma{\cdot}\langle z\rangle$ for some $z\in \Delta {\setminus} \Gamma$. If $KL=\Gamma$, then $\Delta/(K\cap L)\cong (\Gamma/K\times \Gamma/L){\cdot} \langle (Kz,Lz)\rangle\triangleleft_{\,2} \Delta/K\times \Delta/L$.
\end{lemma}

{\bf Proof.} Part (a) follows from a more general Lemma 6.1 and Corollary 6.2 of \cite{J0}, a special case of which (with a different proof) can be found in \cite[Lemma 1]{ACRS}. We therefore prove only the statement (b), based on an argument that covers (a) as well.
\smallskip

To every $h\in \Gamma$ and to every coset of $\Delta$ of the form $(K\cap L)h$ and $(K\cap L)hz$ we assign, respectively, the ordered pair of cosets $(Kh,Lh)$ and $(Khz,Lhz)$ of the group $(\Gamma/K\times \Gamma/L){\cdot}\langle (Kz,Lz)\rangle \triangleleft_{\,2} \Delta/K\times \Delta/L$. Note that $(K\cap L)h = (K\cap L)h'$ for $h,h'\in \Gamma$ is equivalent to $h'h^{-1} = h'z (hz)^{-1}\in K\cap L$ and hence also to both $(Kh,Lh) = (Kh',Lh')$ and $(Khz,Lhz) = (Kh'z,Lh'z)$, so that the above assignment well defines a mapping $\theta:\ \Delta/(K\cap L)\to (\Gamma/K\times \Gamma/L){\cdot}\langle (Kz,Lz)\rangle$. Essentially the same argument shows that $\theta$ is a one-to-one mapping. To show that $\theta$ is a group homomorphism it is sufficient to consider validity of the equality
\begin{equation}\label{eq:parallel}
\theta\left((K\cap L)gz){\cdot}((K\cap L)h\right) = \theta\left((K\cap L)gz\right)\cdot\theta\left((K\cap L)h\right)
\end{equation}
for arbitrary $g,h\in \Gamma$. Since $zh=h'z$ by our assumptions on $\Delta$ and $\Gamma$, the left-hand-side of \eqref{eq:parallel} is equal to $\theta((K\cap L)gh'z) = (Kgh'z,Lgh'z) = (Kgzh,Lgzh) = (Kgz,Lgz)(Kh,Lh)$, which is the right-hand-side of \eqref{eq:parallel}.
\smallskip

It remains to prove that $\theta$ is onto, which we do as follows. Consider an arbitrary element $(Kg,Lh)\in \Gamma/K\times \Gamma/L$. From $KL=\Gamma$ one also has $LK=\Gamma$, and so for every $g,h\in \Gamma$ there exist $u\in K$ and $v\in L$ such that $hg^{-1}=v^{-1}u$, or, equivalently, $ug=vh$. Denoting this last element of $\Gamma$ by $w$, it follows that $Kw=Kug=Kg$ and $Lw=Lvh=Lh$. But then, $\theta((K\cap L)w) = (Kw,Lw)=(Kg,Lh)$. It is clear that for the same $g,h,w$ this argument shows that the element $(K\cap L)wz$ is a pre-image of $(Kgz,Lhz)$ under $\theta$ in the complement of $\Gamma/K\times \Gamma/L$ in the group $(\Gamma/K\times \Gamma/L){\cdot}\langle (Kz,Lz)\rangle$.
\smallskip

This establishes the isomorphism $\Delta/(K\cap L)\cong (\Gamma/K\times \Gamma/L){\cdot}\langle (Kz,Lz)\rangle$ of part (b); the group on the right-hand-side is a subgroup of $\Delta/K\times \Delta/L$ of index $2$. \hfill $\Box$
\medskip

One could give a different argument of part (b) following the method used in the proof of Lemma 1 in \cite{ACRS}, but the one given here has the advantage of containing explicit isomorphisms, which will be used in the next section.

\section{The result}\label{sec:main}

It is a consequence of classification of toroidal orientably-regular maps (see e.g. \cite{JoSi}) that every orientably-regular toroidal map is smoothly covered by infinitely many orientably-regular but chiral toroidal maps. We will therefore consider here only orientably-regular maps of hyperbolic type $\{m,n\}$, that is, such that $1/m+1/n<1/2$.

\begin{theorem}\label{thm:main}
Every orientably-regular map of hyperbolic type is covered by an infinite number of orientably-regular but chiral maps of the same type.
\end{theorem}

{\bf Proof.} For a given hyperbolic type $\{m,n\}$ let $M_1 = {\rm Map}(G_1;x_1,y_1)$ be a chiral, orientably-regular map of type $\{m,n\}$, with $G_1\cong {\rm Aut}^+(M_1)$ isomorphic to a symmetric group $S_r$ of degree $r$, or an alternating group $A_r$ of degree $r$, for some $r\ge 5$; existence of such maps of arbitrary hyperbolic type was established in \cite{CH+}. Let $L_1\triangleleft \Delta=\Delta(m,n,2)$ be the map subgroup of $M_1$, so that $\Delta/ L_1\cong G_1$.
\smallskip

Further, let $M$ be a given orientably-regular map of an arbitrary hyperbolic type $\{m,n\}$, as in the statement. By residual finiteness of triangle groups (see, for example, \cite[Section 3.4]{Si-surv}) applied to the map subgroup of $M$, or by the Macbeath trick (as used, for example, in \cite{CH+}), there exist infinitely many orientably-regular maps that smoothly cover $M$. Let $M_2 = {\rm Map}(G_2;x_2,y_2)$ be such a smooth cover of $M$; because of their infinite inventory we may assume that the group $G_2$ is not isomorphic to $G_1$. Also, let $L_2 \triangleleft \Delta$ be the map subgroup of $M_2$, so that $\Delta/L_2\cong G_2$.
\smallskip

We will consider the parallel product $M_1||M_2= {\rm Map}(G;x,y)$ of our `seed' chiral map and a smooth cover of the given map. The map $M_1||M_2$ is orientably-regular, of type $\{m,n\}$, and its map subgroup is $L_1\cap L_2\, \triangleleft\, \Delta$, so that $G\cong \Delta/(L_1\cap L_2)$. Our aim is to show that $M_1||M_2$ is chiral, which amounts to proving that it admits no orientation-reversing automorphism. We saw that this is equivalent to non-existence of an automorphism of the group $G={\rm Aut}^+ (M_1||M_2)$ inverting both generators $x,y$ of $G$. We will divide the proof into two parts.
\smallskip

{\bf Part 1.} Let $L_1L_2= \Delta = \Delta(X,Y)$. By Lemma \ref{lem:parallel}, the group $\Delta/(L_1\cap L_2)$ is isomorphic to the direct product $\Delta/L_1\times \Delta/L_2$, via an isomorphism given by $(L_1\cap L_2)X\mapsto (L_1X,L_2X)$ and $(L_1\cap L_2)Y\mapsto (L_1Y,L_2Y)$. In terms of the generators $x,y$ of $G\cong \Delta/(L_1\cap L_2)$ and $x_i,y_i$ of $G_i\cong \Delta/L_i$ for $i=1,2$, this means that the group $G$ is isomorphic to $G_1\times G_2$ through an isomorphism given by $x\mapsto (x_1,x_2)$ and $y\mapsto (y_1, y_2)$.
\smallskip

To obtain more information about the structure of the automorphism group ${\rm Aut}(G)$ of the group $G\cong G_1\times G_2$, we will invoke Corollaries 3.8 and 3.10 of \cite{BCM}, both of which assume that $G_1$ and $G_2$ have no common direct factor. This is satisfied in our setting, because $G_1$ is symmetric or alternating (of degree at least $5$) and $G_2\not\cong G_1$; if the two groups contained a common direct factor, then $G_1$ would have to be a direct factor of $G_2$, which is impossible since $G_2$ is $2$-generated.
\smallskip

If $G_1\cong A_r$ for some $r\ge 5$, which is a simple group, one may apply \cite[Corollary 3.8]{BCM} to conclude that ${\rm Aut}(G)\cong {\rm Aut}(G_1)\times {\rm Aut}(G_2)$, so that in this case the only automorphisms of $G\cong G_1\times G_2$ are ordered pairs of automorphisms of $G_1$ and $G_2$. If $G_1\cong S_r$ for some $r\ge 5$, then, by a special case of \cite[Corollary 3.10]{BCM} applied to the simple subgroup $A_r\triangleleft_{\,2} S_r$, it follows that every automorphism of the product $S_r\times G_2$ has the form $(\sigma,\gamma;c^j)$ for some $\sigma\in {\rm Aut} (S_r)$, $\gamma\in {\rm Aut}(G_2)$, $c$ a central involution in $G_2$, and $j\in \{0,1\}$. The action of such an automorphism on pairs $(\pi,g)\in S_r\times G_2$ is, by \cite{BCM}, given by $(\pi,g)(\sigma,\gamma;c^j)= (\pi\sigma,g\gamma c^{j_\pi})$, where $j_\pi = 0$ if $j=0$ (for every $\pi\in S_r$) and also if $j=1$ but only for every $\pi \in A_r$, while $j_\pi=1$ if $j=1$ and $\pi$ is an arbitrary element of $S_r{\setminus}A_r$. The important conclusion is that, in both possibilities for $G_1$, the projection of every automorphism of $G\cong G_1\times G_2$ onto $G_1$ is an automorphism of $G_1$.
\smallskip

Suppose now that $\tau$ is an orientation-reversing automorphism of $M=M_1||M_2$, that is, an automorphism of the group $G$ such that $x\tau=x^{-1}$ and $y\tau = y^{-1}$. Since $G\cong G_1\times G_2$ by an isomorphism sending $x$ to $(x_1,x_2)$ and $y$ to $(y_1,y_2)$, we may consider $\tau$ to be an automorphism of $G_1\times G_2$ sending $(x_1,x_2)$ to $(x_1^{-1},x_2^{-1})$ and $(y_1,y_2)$ to $(y_1^{-1},y_2^{-1}) $. This together with findings of the previous paragraph shows that $\tau$ projects to an automorphism of $G_1$ that inverts the generators $x_1$ and $y_1$ of $G_1$. But this contradicts the fact that the map $M_1= {\rm Map}(G_1; x_1,y_1)$ we started with was chiral.
\medskip

{\bf Part II.} Let now $L_1L_2$ be a proper subgroup of $\Delta$. From $G_1\not\cong G_2$ one has $L_1\ne L_2$,  implying that the group $\Gamma=L_1L_2$ is a normal subgroup of $\Delta$ properly containing $L_1$. The normal series $L_1\,\triangleleft\, \Gamma \,\triangleleft\, \Delta$ corresponds to the series $1\triangleleft\, \Gamma/L_1 \triangleleft\, \Delta/L_1\cong G_1$, with all the normal inclusions proper. The only possibility then is that $G_1\cong S_r$ and $H_1:=\Gamma/L_1\cong A_r$ by the Correspondence Theorem, which also implies that $\Gamma$ is a (normal) subgroup {\em of index} $2$ in $\Delta$. By Section \ref{sec:special}, a normal subgroup of $\Delta= \Delta(m,n)$ of index $2$ exists if and only if one of $m$, $n$ is even; taking duality into the account we will assume that $n$ is even.
\smallskip

By the analysis of subgroups of $\Delta$ of index $2$ in section \ref{sec:special}, we may identify the group $\Gamma= L_1L_2$ with one of $\Gamma_A =\Gamma(X,Y^2)$, $\Gamma_B =\Gamma(X^2,Y^2)$, according to the type (A or B) of $\Gamma$. In the terminology of section \ref{sec:special}, not only both $M_1=(G_1;x_1,y_1)$ and $M_2=(G_2;x_2,y_2)$, but also their parallel product $M_1||M_2={\rm Map}(G;x,y)$, are {\em $\Gamma$-exceptional} (the last one because its map subgroup $L_1\cap L_2$ is also contained in $\Gamma$). It follows that the normal subgroups $H\cong \Gamma/(L_1\cap L_2)$, $H_1\cong \Gamma/L_1$ and $H_2\cong \Gamma/L_2$, of index $2$ in $G$, $G_1$ and $G_2$, respectively, admit a presentation of the form \eqref{eq:H-Y2} if $\Gamma=\Gamma_A$ (in terms of the corresponding generators $u=yxy^{-1}$, $v=y^2$ and $u_i=y_ix_iy_i^{-1}$, $v_i=y_i^2$, $i=1,2$), or in the form \eqref{eq:H-X2Y2} if $\Gamma=\Gamma_B$ (in terms of $u=xy^{-1}$, $v=y^2$, $w=yx$ and $u_i=x_iy_i^{-1}$, $v_i=y_i^2$, $w_i=y_ix_i$, $i=1,2$).
\smallskip

Observe that, for $i=1,2$, the groups $G$, $G_i$ and $H$, $H_i$ are related by  $G=H\langle y\rangle$ and $G_i=H_i\langle y_i\rangle$, with $y^2\in H$, $y_i^2\in H_i$ for both types of $\Gamma$. In the product, the element  $y\in G{\setminus}H$ acts on $H$ by $y^{-1}uy=(uv)^{-1}$ and $y^{-1}vy=v$ for type A, and by $y^{-1}uy=v^{-1}w$, $y^{-1}vy=v$ and $y^{-1}wy=uv$ for type B, with completely analogous actions on $H_1$ and $H_2$. Now, Lemma \ref{lem:parallel}, applied with $y$ and $y_i$ in place of $z$ to quotients of $\Gamma=L_1L_2$ and $\Delta\cong \Gamma\langle Y\rangle$ (with an analogous action of $Y$ on $\Gamma$ as above), implies that $H\cong H_1\times H_2$ and $G=H\langle y\rangle \cong (H_1\times H_2)\langle (y_1,y_2)\rangle\, \triangleleft_{\,2}\, H_1\langle y_1 \rangle \times H_2\langle y_2\rangle = G_1\times G_2$. Moreover, the proof of Lemma \ref{lem:parallel} shows the following, in terms of cosets of $\Gamma$ (and hence also of $\Delta$):

\begin{quote}\label{quo:Delta}
(*) {\sl With $\Delta=\Gamma\langle Y\rangle$, an isomorphism $\Delta/(L_1\cap L_2) \to (\Gamma/L_1 \times \Gamma/L_2) \langle (L_1Y,L_2Y)\rangle$ is given, for every $Z\in \Gamma$, by sending the cosets $(L_1\cap L_2)Z\in H=\Gamma/ (L_1\cap L_2)$ and $(L_1\cap L_2)ZY\in G{\setminus}H$ to the ordered pairs $(L_1Z,L_2Z)\in H_1\times H_2$ and $(L_1ZY,L_2ZY)\in (H_1\times H_2)\langle (y_1,y_2)\rangle \setminus (H_1\times H_2)$.     }
\end{quote}

\noindent A translation of (*) to the groups $H\, \triangleleft_{\,2}\, G$ and $H_i\, \triangleleft_{\,2}\, G_i$ ($i=1,2$) reads as follows:

\begin{quote}\label{quo:G}
(**) {\sl Letting $z=(L_1\cap L_2)Z\in H$ and $z_i = L_iZ\in H_i=\Gamma/L_i$, the isomorphism $H\langle y\rangle \cong (H_1\times H_2)\langle (y_1,y_2)\rangle$ from the proof of part (b) in Lemma \ref{lem:parallel} has the form $z\in H \mapsto (z_1,z_2)\in H_1\times H_2$ and $zy\in G{\setminus}H \mapsto (z_1y,z_2y)$, where the last pair is an element of the group $(H_1\times H_2)\langle (y_1,y_2)\setminus (H_1\times H_2)$.}
\end{quote}

\noindent The way $z$ and $z_i$ depend on $Z$ and $Z_i$ ($i=1,2$) has numerous consequences. For example, if $\Gamma$ is of type B, then $z=u\in H$ (for $Z=U\in \Gamma$) gives $zy=x\in G{\setminus}H$ (corresponding to $ZY=X \in \Delta{\setminus}\Gamma$); similarly, $z_i=u_i\in H_i$ gives $z_iy_i=x_i\in G_i{\setminus}H_i$ for $i=1,2$. The isomorphism in (**) then sends $u\in H$ to $(u_1,u_2)\in H_1\times H_2$, and $x\in G{\setminus}H$ to $(x_1,x_2)$, an element of $(H_1\times H_2)\langle (y_1,y_2)\setminus (H_1\times H_2)$. Informally, (**) means that all relations implied by dependence of $U,V,W\in \Gamma$ on $X,Y\in \Delta$ project onto dependence relations of $u,v,w\in H$ on $x,y\in G$, including their indexed versions for $i=1,2$ in $H_i$ and $G_i$.
\smallskip

We are finally prepared to take up the main step. Suppose that the map $M_1||M_2$ is reflexible, and let $\tau$ be and automorphism of the group $G=\langle x,y\rangle \cong \Delta/(L_1\cap L_2)$ inverting both $x$ and $y$. Since $u\tau = v^{-1}u^{-1}v$ if $\Gamma=\Gamma_A$, $u\tau = wv$ and $w\tau = (uv)^{-1}$ if $\Gamma=\Gamma_B$, and $v\tau = v^{-1}$ in both cases, the restriction $\tau_H$ of $\tau$ to $H$ is an automorphism of $H$. We will carry on with details only for $\Gamma$ of type B; the arguments for type A are analogous.
\smallskip

As seen earlier in this part of the proof, in the situation when $\Gamma=\Gamma_B=L_1L_2$ one has $H\cong H_1\times H_2$ by Lemma \ref{lem:parallel} (a). But this time, $H_1=\Gamma/L_1\cong A_r$ for $r\ge 5$, which is a simple group. Applying Corollary 3.8 of \cite{BCM} then gives ${\rm Aut}(H)\cong {\rm Aut}(H_1)\times {\rm Aut}(H_2)$, so that the only automorphisms of $H\cong H_1\times H_2$ are ordered pairs of automorphisms of $H_1$ and $H_2$. In particular, for type B it follows that the automorphism $\tau_H$ sending the ordered triple $(u,v,w)$ to $(wv,v^{-1},(uv)^{-1})$ induces an automorphism $\tau_1$ of $H_1$, mapping the triple $(u_1,v_1,w_1)$ to $(w_1v_1,v_1^{-1},(u_1v_1)^{-1})$.
\smallskip

Recalling the consequences of our notation introduced in (**), the generators $u_1,v_1,w_1$ of $H_1$ are tied to the generators $x_1,y_1$ of $H_1\langle y_1\rangle = G_1$ by $u_1=x_1y_1^{-1}$, $v_1=y_1^2$ and $w_1=y_1x_1$. But now, part (b) of Lemma \ref{lem:HtoG} applies, with the conclusion that the map $M_1$ is reflexible, giving a contradiction and completing the proof. \hfill $\Box$
\bigskip

The reader may have noticed that the proof of Theorem \ref{thm:main} would be much shorter if one had an inventory of chiral orientably-regular maps of an arbitrary hyperbolic type but with an {\em alternating} automorphism group. In such a case all one would need for the parallel product construction is part (a) of Lemma \ref{lem:parallel} in combination with Corollary 3.8 of \cite{BCM}. This is actually the case for all hyperbolic types $\{m,n\}$ in which both entries are odd: if $M={\rm Map}(G;x,y)$ is a map of such a type with $G$ isomorphic to a symmetric or an alternating group, then $G$ is necessarily alternating because of even parity of the generating permutations $x$ and $y$ (both of odd order). Thus, we actually have a very short proof of Theorem \ref{thm:main} for asymptotically a quarter of all hyperbolic types.
\smallskip

To the best of our knowledge, no construction of chiral orientably-regular maps of an arbitrary given type with alternating automorphism group was available at the time of submission of this note. Such a construction would have even more striking applications if it was available for {\em hypermaps}, which can be identified with finite quotients of $(m,n,\ell)$-triangle groups $\Delta(m,n,\ell) = \langle X,\,Y\ |\ X^m,\,Y^n,\,(XY)^\ell\rangle$. In such a case one would have a very short proof of the full version of the result of Jones \cite{J1}, because the shortcut alluded to in the previous paragraph would apply to hypermaps as well.

\bigskip

{\bf Acknowledgment.} The second author acknowledges support from APVV Research Grants 19-0308 and 22-0005, and from VEGA Research Grants 1/0567/22 and 1/0069/23.

\end{document}